\documentclass[preprint,12pt]{elsarticle}

\usepackage{hyperref,amsmath,amsfonts,mathrsfs}

\usepackage{graphicx}
\usepackage{amsfonts}
\usepackage{amsbsy}
\usepackage{amssymb}
\usepackage{amsmath,tikz,pgfplots}
\usetikzlibrary{decorations.pathmorphing}
\usepackage{tikz}
\numberwithin{equation}{section}
\usetikzlibrary{arrows}
\usepackage{graphics}
\usepackage{graphicx}
\usepackage{amsfonts}
\usepackage{amsbsy}
\usepackage{amssymb}
\usepackage{amsmath,tikz,pgfplots}
\usetikzlibrary{decorations.pathmorphing}
\usepackage{graphics}
\def\bpsp{\begin{pspicture}}
\def\epsp{\end{pspicture}}

\newtheorem{theorem}{Theorem}[section]
\newtheorem{remark}[theorem]{Remark}
\newtheorem{example}[theorem]{Example}
\newtheorem{lemma}[theorem]{Lemma}
\newtheorem{corollary}[theorem]{Corollary}
\newtheorem{definition}[theorem]{Definition}
\newtheorem{proposition}[theorem]{Proposition}
\newtheorem{note}{Note}
\newtheorem{case}{Case}
\newtheorem{conjecture}{Conjecture}
\newtheorem{question}{Question}

\newcommand{\bea}{\begin{eqnarray}}
\newcommand{\eea}{\end{eqnarray}}
\newcommand{\beq}{\begin{eqnarray*}}
\newcommand{\eeq}{\end{eqnarray*}}

\def\m4{\mbox{\rm ~(mod $4$)}}

\def \bd{\begin{definition}}
\def \ed{\end{definition}}
\def \bqu{\begin{question}}
\def \equ{\end{question}}
\def \bcc{\begin{conjecture}}
\def \ecc{\end{conjecture}}
\def \bt{\begin{theorem}}
\def \et{\end{theorem}}
\def \bl{\begin{lemma}}
\def \el{\end{lemma}}
\def \bc{\begin{corollary}}
\def \ec{\end{corollary}}
\def \be{\begin{equation}}
\def \ee{\end{equation}}
\def \ben{\begin{enumerate}}
\def \een{\end{enumerate}}
\def \ba{\begin{array}}
\def \ea{\end{array}}
\def \bp{\begin{proposition}}
\def \ep{\end{proposition}}
\def \bx{\begin{example}}
\def \ex{\end{example}}
\def \br{\begin{remark}}
\def \er{\end{remark}}
\def \bdsc{\begin{description}}
\def \edsc{\end{description}}

\def \bn{\begin{case}}
\def \en{\end{case}}
\def \bnt{\begin{note}}
\def \ent{\end{note}}
\def\1{1\!\!1}

\def\mm2{\mbox{\rm ~(mod $2$)}}
\def\m4{\mbox{\rm ~(mod $4$)}}

\def\qed{\nolinebreak\hfill\rule{.2cm}{.2cm}\par\addvspace{.5cm}}

\def\m{\mu}

\def\1{\textbf{1}}
\def\0{\textbf{0}}

\journal{ABC}

\begin{document}

\begin{frontmatter}



\title{Bounds for the trace norm of $A_{\alpha}$ matrix of digraphs}
\author{Mushtaq A. Bhat$^{a}$}
\author{Peer Abdul Manan$^{b}$}

\address{Department of  Mathematics,  National Institute of Technology, Srinagar-190006, India}
\address {$^{a}$mushtaqab@nitsri.ac.in;~~~$^{b}$\text{mananab214@gmail.com}}

\begin{abstract}
 Let $D$ be a digraph of order $n$ with adjacency matrix $A(D)$. For $\alpha\in[0,1)$, the $A_{\alpha}$ matrix of $D$ is defined as $A_{\alpha}(D)=\alpha {\Delta}^{+}(D)+(1-\alpha)A(D)$, where ${\Delta}^{+}(D)=\mbox{diag}~(d_1^{+},d_2^{+},\dots,d_n^{+})$ is the diagonal matrix of vertex outdegrees of $D$. Let $\sigma_{1\alpha}(D),\sigma_{2\alpha}(D),\dots,\sigma_{n\alpha}(D)$ be the singular values of $A_{\alpha}(D)$. Then the trace norm of $A_{\alpha}(D)$, which we call $\alpha$ trace norm of $D$, is defined as $\|A_{\alpha}(D)\|_*=\sum_{i=1}^{n}\sigma_{i\alpha}(D)$. In this paper, we find the singular values of some basic digraphs and characterize the digraphs $D$ with $\mbox{Rank}~(A_{\alpha}(D))=1$. As an application of these results, we obtain a lower bound for the trace norm of $A_{\alpha}$ matrix of digraphs and determine the extremal digraphs. In particular, we determine the oriented trees for which the trace norm of $A_{\alpha}$ matrix attains minimum. We obtain a lower bound for the $\alpha$ spectral norm $\sigma_{1\alpha}(D)$ of digraphs and characterize the extremal digraphs. As an application of this result, we obtain an upper bound for the $\alpha$ trace norm of digraphs and characterize the extremal digraphs.
\end{abstract}

\vskip 0.2 true cm


\end{frontmatter}

\section{\bf Introduction}
A directed graph (or briefly digraph) $D$ consists of two sets $\mathcal{V}$ and $\mathcal{A}$, where $\mathcal{V}$ is a non-empty finite set whose elements are called vertices and $\mathcal{A}$ is a set of ordered pairs of elements of $\mathcal{V}$ and is known as a set of arcs. We assume our digraphs are simple i.e., there are no loops and parallel arcs. A graph $G$ can be identified with a symmetric digraph $\overleftrightarrow {G}$ obtained by replacing each edge $e$ of $G$ by a pair of symmetric arcs. We call a digraph to be asymmetric if it has no pair of symmetric arcs. An asymmetric digraph is also known as an oriented graph. In a digraph $D$, an arc from a vertex $u$ to $v$ is denoted by $(u,v)$. In this case, we say $u$ is the tail and $v$ is the head of arc $(u,v)$. The set of vertices $\{w\in \mathcal{V}: (u,w)\in \mathcal{A}\}$ is called the outer neighbour of $u$ and we denote it by $N^{+}(u)$. The set of vertices $\{w\in \mathcal{V}: (w,u)\in \mathcal{A}\}$ is called the inner neighbour of $u$ and we denote this by $N^{-}(u)$. The cardinality of the set $N^{+}(u)$ is called outdegree of $u$ and we denote it by $d_u^+$. Therefore, $d_u^+=|N^{+}(u)|$. The cardinality of the set $N^{-}(u)$ is called indegree of $u$ and we denote it by $d_u^-$. Therefore, $d_u^-=|N^{-}(u)|$. \\
\indent Let $D$ be a digraph with vertex set $\mathcal{V}=\{v_1,v_2,\dots,v_n\}$. Then the adjacency matrix $A(D)=(a_{ij})$ of $D$ is a square matrix of order $n$ with $a_{ij}=1$ if there is an arc from vertex  $v_i$ to vertex $v_j$ and zero, otherwise. Let ${\Delta}^+={\Delta}^+(D)=\mbox{diag}~(d^+_1,d^+_2,d^+_3,\dots,d^+_n)$, where $d^+_i=d^+_{v_i}$, be the the outdegree matrix of $D$. Then the Laplacian and signless Laplacian matrices of $D$ are respectively defined as $L(D)={\Delta}^{+}+A(D)$ and $Q(D)={\Delta}^{+}-A(D)$. For $\alpha\in[0, 1]$, Nikiforov \cite{n3} defined the alpha adjacency matrix $A_{\alpha}(G)$ of a graph $G$ as a common extension of adjacency matrix and signless Laplacian matrix as $A_{\alpha}(G)=\alpha D(G)+(1-\alpha)A(G)$, where $D(G)$ is diagonal matrix of vertex degrees of $D$. The spectral properties of alpha matrix of graphs are well studied, for example see \cite{lds,llx} and references cited therein.  Liu et al. \cite{lwcl}, defined the $A_{\alpha}$ matrix for directed graphs and studied $\alpha$-spectral radius of digraphs. For $\alpha\in [0, 1)$, the alpha matrix of a digraph is defined as $A_{\alpha}(D)={\Delta}^+(D)+(1-\alpha)A(D)$. For $\alpha=0$, we see $A_0=A$ and for $\alpha=\frac{1}{2}$, $A_{\frac{1}{2}}=\frac{1}{2}Q(D)$. It is clear that $A_{\alpha}(D)$ is a common extension of adjacency matrix $A=A(D)$ and signless Laplacian matrix $Q(D)$ of a digraph $D$. Very recently, Yang et al. \cite{ybw} studied the spectral moments of $A_{\alpha}$ adjacency matrix of a digraph.\\
Let $G$ be an undirected graph of order $n$ and with eigenvalues $\lambda_1\ge\lambda_2\ge\lambda_3\ge\dots\ge\lambda_n$. Then the energy of $G$ is defined as $E(G)=\sum\limits_{k=1}^{n}|\lambda_k|$. This concept was given by Gutman (1978). For details related to graph energy see \cite{lsg}. The concept of energy was extended to digraphs by Pe\~{n}a and Rada \cite{pr} and they defined the energy of a digraph $D$ as $E(D)=\sum_{i=1}^{n}|\Re z_i|$, where $z_1,z_2,\dots,z_n$ are eigenvalues of $D$, possibly complex and $\Re z_i$ denote the real part of complex number $z_i$. Recall that the trace norm of a complex matrix $B\in M_n(\mathbb{C})$ is defined as 
\begin{equation*}
\|{B} \|_* =\sum\limits_{i=1}^{n}\sigma_i(B),
\end{equation*}
where $\sigma_1(B)\ge\sigma_2(B)\ge \sigma_3(B)\ge \dots\ge \sigma_n(B)$ are the singular values of $B$ i.e., the positive square roots of eigenvalues of $BB^{*}$.\\
Throughout paper, we call singular values of $A_{\alpha}(D)$ as the $\alpha$ singular values of $D$ and trace norm of $A_{\alpha}(D)$ as the trace norm of $D$ and we will denote the trace norm of $A_{\alpha}(D)$ by $\|D_{\alpha}\|_{*}$. If a digraph $D$ has $k$ distinct singular values $\sigma_{1\alpha}, \sigma_{2\alpha}, \dots, \sigma_{k\alpha}$, with their respective multiplicities $m_1, m_2,\dots,m_k$, then we write the set of singular values as $$\{{\sigma^{[m_1]}_{1\alpha}}, {\sigma^{[m_2]}_{2\alpha}},\dots,{\sigma^{[m_k]}_{k\alpha}}\}.$$
If $B=A(G)$, the adjacency matrix of graph $G$, then $\sigma_i(B)=|\lambda_i(G)|$ and so trace norm coincides with the energy of a graph. Trace norm of a matrix is also known as Nikiforov energy of a matrix \cite{n1}. So, graph energy extends to digraphs via trace norm as well.\\
Kharaghani and Tayfeh-Rezaie \cite{kr} obtained upper bounds on the trace norm of $(0,1)$-matrices. Agudelo and Rada \cite{ar} obtained lower bounds for the trace norm (Nikiforov's energy) of adjacency matrices of digraphs. Agudelo, Pe\~{n}a and Rada \cite{apr} determined trees attaining minimum and maximum trace norm. Monsalve and Rada \cite{mr} determined oriented  bipartite graphs with minimum trace norm. Agudelo, Rada and Rivera \cite{arr} obtained upper bound for the trace norm of the Laplacian matrix of a digraph in terms of number of vertices $n$, number of arcs $a$ and out degrees of a digraph. For more about trace norm see \cite{jgmr,n2,n3,psg}.\\
The rest of the paper is organized as follows.\\
In section $2$, we obtain the singular values of the $A_{\alpha}$ matrix of a directed path $\overrightarrow{P_n}$, with arc set $\mathcal{A}(\overrightarrow{P_n})=\{(v_1,v_2),(v_2,v_3),\dots,(v_{n-1}, v_n)\}$, directed cycle $\overrightarrow{C_n}$, with arc set $\mathcal{A}(\overrightarrow{C_n})=\{(v_1,v_2),(v_2,v_3),\dots,(v_{n-1}, v_n),(v_n,v_1)\}$,  and ${\overrightarrow{K}}_{r,s}$, where ${\overrightarrow{K}}_{r,s}$ denote oriented complete bipartite graph with partite sets $\{u_1,u_2,\dots,u_r\}$, $\{v_1,v_2,\dots,v_s\}$ and all arcs of the form $(u_i,v_j)$, where $i=1,2,\dots,r$ and $j=1,2,\dots,s$. We characterize digraphs with $Rank(A_{\alpha}(D))=1$. Using these results, we present a lower bound for the trace norm of $A_{\alpha}(D)$ of a digraph in terms of $n,a$ and sum of squares of outdegrees of $D$ and we characterize the digraphs attaining the lower bound. As a consequence, we determine oriented trees with minimum trace norm for $A_{\alpha}$ matrix.\\
In section $3$, we obtain upper bounds for trace norm of $A_{\alpha}$ matrix of digraphs and characterize the extremal digraphs. 
\section{\bf Lower bounds for trace norm of $A_{\alpha}$ matrix of digraphs}\label{sec2}
We first compute $\alpha$-singular values of $\overrightarrow{P_n},\overrightarrow{C_n}$ and $\overrightarrow{K}_{r,s}$, which will be used in our main results.
\begin{lemma} \label{2.1}
 If $\overrightarrow{P_n}$ denote a directed path of order $n$, then the $\alpha$-singular values of $\overrightarrow{P}_n$ are 
$$0^{[1]} ~~~\mbox{and} ~~\sqrt{2\alpha^2-2\alpha+1+2\alpha(1-\alpha)\cos{\frac{j\pi}{n}}},$$ where $j=1,2,\ldots,n-1$.
\end{lemma}
 {\bf Proof.} We have $A_{\alpha}({\overrightarrow{P_{n}}})=\alpha {\Delta^{+}}({\overrightarrow{P_{n}}})+(1-\alpha)A({\overrightarrow{P_{n}}})$,\\ so that 
 \begin{align*} 
 A_{\alpha}{A_{\alpha}}^{T}&=(\alpha {\Delta^{+}}+(1-\alpha)A)(\alpha {\Delta^{+}}+(1-\alpha)A^{T}) \\
 &=(\alpha {\Delta^{+}}+(1-\alpha)A)(\alpha {\Delta^{+}}+(1-\alpha)A^{T}) \\
  &=\alpha^2 ({\Delta^{+}})^2+\alpha(1-\alpha){\Delta^{+}}A^{T}+\alpha(1-\alpha)A{\Delta^{+}}+(1-\alpha)^2AA^{T} \\
  &=\alpha^2 ({\Delta^{+}})^2+\alpha(1-\alpha)(A{\Delta^{+}}+{\Delta^{+}}A^{T})+(1-\alpha)^2AA^{T}\\ 
    &=(\alpha^2+(1-\alpha)^2)\begin{bmatrix} 
I_{n-1} & {\bf{0}}_{n-1 \times 1}\\
{\bf{0}}_{1\times n-1}& 0_{1 \times 1}
\end{bmatrix}
+ \alpha(1-\alpha)\begin{bmatrix} 
A({\overleftrightarrow{P_{n-1}}}) & {\bf{0}}_{n-1 \times 1}\\
{\bf{0}}_{1\times n-1}& 0_{1 \times 1}
\end{bmatrix}
\end{align*}
As the eigenvalues of path $P_{n-1}$ are $2\cos \frac{\pi j}{n}$, where $j=1,2,\dots,n-1$. Consequently, the eigenvalues of $A_{\alpha}{A_{\alpha}}^{T}$ are\\
$$0^{[1]} ~~~\mbox{and} ~~2\alpha^2-2\alpha+1+2\alpha(1-\alpha)\cos{\frac{j\pi}{n}},$$ where $j=1,2,\ldots,n-1$.\\
The singular values of $A_{\alpha}$ are\\
$$0^{[1]} ~~~\mbox{and} ~~\sqrt{2\alpha^2-2\alpha+1+2\alpha(1-\alpha)\cos{\frac{j\pi}{n}}},$$ where $j=1,2,\ldots,n-1$. \qed
\begin{lemma} \label{2.2}
 If $\overrightarrow{C_{n}}$ is a directed cycle of order $n$, then the $\alpha$-singular values of $\overrightarrow{C_{n}}$ are 
$$
\sqrt{(2\alpha^{2}-2\alpha+1)+2\alpha(1-\alpha)\cos\frac{2j\pi}{n}},$$
 where $j=0,1,2,3,\dots,n-1.$
 \end{lemma}
{\bf Proof.}
We have $A_{\alpha}({\overrightarrow{C_{n}}})=\alpha {\Delta^{+}}({\overrightarrow{C_{n}}})+(1-\alpha)A({\overrightarrow{C_{n}}})$, so that
\begin{align*}
A_{\alpha}{A_{\alpha}}^{T}&=(\alpha {\Delta^{+}}+(1-\alpha)A)(\alpha {\Delta^{+}}+(1-\alpha)A^{T})\\
&=(\alpha I_{n}+(1-\alpha)A)(\alpha I_{n}+(1-\alpha)A^{T})\\
&=\alpha^{2}I_{n}+\alpha(1-\alpha)A^{T}+\alpha(1-\alpha)A+(1-\alpha)^2AA^{T}\\
&=[\alpha^2+(1-\alpha)^2]I_n+\alpha(1-\alpha)(A+A^{-1}),\\
 \end{align*}
by using the fact that $AA^{T}=I_n$, for adjacency matrix $A$ of directed cycle ${\overrightarrow{C_{n}}}$. We see that the eigenvalues of $A_{\alpha}{A_{\alpha}}^{T}$ are \\
$(2\alpha^{2}-2\alpha+1)+\alpha(1-\alpha)(\omega^{j}+\frac{1}{\omega^{j}})$,
where $\omega^{n}=1$ and $j=0,1,2,\dots,n-1.$\\
 $\mbox{or} ~~~\sigma_{j}^{2}=(2\alpha^{2}-2\alpha+1)+\alpha(1-\alpha)\cos\frac{2j\pi}{n}$,
 where $j=0,1,2,\dots,n-1.$\\
$ \mbox{or}~~~ \sigma_{j}=\sqrt{(2\alpha^{2}-2\alpha+1)+\alpha(1-\alpha)\cos\frac{2j\pi}{n}}$,
  where $j=0,1,2,\dots,n-1.$ \qed

\begin{lemma}\label{2.3}
If ${\overrightarrow{K}}_{r,s}$ be a complete bipartite digraph with partite sets $X=\{u_1,u_2,\dots,u_r\}$ and $Y=\{v_1,v_2,...,v_s\}$ and arcs $(u_i,v_j)$  where $i=1,2,\dots,r$ and $j=1,2,\dots,s$, then the $\alpha$-singular values of  ${\overrightarrow{K}}_{r,s}$ are
$$0^{[s]},~~(\alpha s)^{[r-1]}~~~~and~~~~~\sqrt{\alpha^2s^2+(1-\alpha)^2sr}$$
\end{lemma}
{\bf Proof.} We have
\[ A_{\alpha}({\overrightarrow{K}}_{r,s})=\begin{bmatrix} 
\alpha sI_{r} & (1-\alpha)J_{r \times s}\\
{\bf{0}}_{s\times r}& {\bf{0}}_{s \times s}
\end{bmatrix}
\]
 so that \[ A_{\alpha}A_{\alpha}^{T}=\begin{bmatrix} 
\alpha^2s^2I_{r}+(1-\alpha)^2sJ_{r\times r} & {\bf{0}}_{r \times s}\\
{\bf{0}}_{s\times r}& {\bf{0}}_{s \times s}
\end{bmatrix}
\]
The eigenvalues of $A_{\alpha}{A_{\alpha}}^{T}$ are 
$$0^{[s]},~~~(\alpha^2s^2)^{[r-1]} ~~~~~\mbox{and}~~~~~ \alpha^2s^2+(1-\alpha)^2sr$$
Consequently, the singular values of $A_{\alpha}$ are
$$0^{[s]},~~(\alpha s)^{[r-1]}~~~\mbox{and}~~~\sqrt{\alpha^2s^2+(1-\alpha)^2sr}.$$\qed 

\begin{remark}\label{2.4} For $\alpha\in(0,1)$, $[A_{\alpha}(D)]^{T}$ need not be same as $A_{\alpha}(D^{T})$, where $D^{T}$ denotes transpose (or converse) of digraph $D$. So a digraph and its transpose need not have same $\alpha$-singular values.
\end{remark}
Recall that the rank of a matrix $B\in M_{m\times n}(\mathbb{C})$ is the number of nonzero singular values of $B$. Rank of a matrix is also defined as the maximum number of its linearly independent rows or columns. Recently rank of digraphs has been studied in \cite{amr,zxw}. The authors define rank preserving operations and characterize oriented graphs with rank $1,2,3$. We characterize digraphs $D$ such $\mbox{Rank}(A_{\alpha}(D))=1$ and this result will be used to discuss equality case in our main results. We give a different proof here, which will serve as an alternative proof to corresponding rank characterization for adjacency matrix of digraphs.
\begin{theorem}\label{2.5}
If $D$ is a digraph of order $n\ge 2$. Then $\mbox{Rank}(A_{\alpha}(D))=1$ if and only if $\alpha=0$ and $D={\overrightarrow{K}}_{r,s}+$~Possibly some isolated vertices or $\alpha=\frac{1}{2}$ and $D=\overleftrightarrow{K_2}$.
\end{theorem}
{\bf Proof.} We prove this result when $D$ has no isolated vertices, since adding isolated vertices does not change the rank. For $n=2$, it is easy to see that $\mbox{Rank}(A_{\alpha}(D))=1$, if and only if $\alpha=\frac{1}{2}$ and $D=\overleftrightarrow{K_2}$ or $\alpha=0$ and $D=\overrightarrow{P_2}$. Assume $n\ge 3$ and $\mbox{Rank}(A_{\alpha}(D))=1$. We first claim that $D$ cannot have two consecutive arcs of the form $(v_1,v_2)$ and $(v_2,v_3)$, where $X=\{v_1,v_2,v_3\}\subseteq \mathcal{V}(D)$. For if there are two consecutive arcs $(v_1,v_2)$ and $(v_2,v_3)$. Then $D$ contains one of the eleven digraphs $D_i$, where $i=1,2,\dots,11$ shown in Fig. $1$ as induced subdigraphs on vertex subset $X$.\\
If $D$ contains $D_1,D_2 \mbox~{or}~D_3$, then since there can be arcs from $X$ to $\mathcal{V}-X$, the principal submatrix of $A_{\alpha}(D)$ corresponding to vertices in $X$ is of the form\\

\[
{P}=
\begin{bmatrix} 
\alpha d_1^+ & 1-\alpha & 0\\
0 & \alpha d_2^+ & 1-\alpha \\
0 & 0 & \alpha d_3^+
\end{bmatrix}
\]

with $d_1^+\ge 1, d_2^+\ge 1, d_3^+\ge 0 $
and $2 \le d_1^+ +d_2^+ +d_3^+ \le a$\\
if $\alpha \ne 0$, then $Rank(P)=2$ or $3$ according as $d_3^+=0$ or not.\\
Also, for $\alpha =0$,  $\mbox{Rank}(P)=2$. Consequently, $\mbox{Rank}(A_{\alpha}(D))$ is atleast $2$ in this case, a contradiction.\\
If $D$ contains one of induced subdigraphs $D_4,D_5,D_6 \mbox~{or}~D_7$, then the principal submatrix $Q$ of $A_{\alpha}(D)$ corresponding to vertices in $X$ has the form\\

\[
{Q}=
\begin{bmatrix} 
\alpha d_1^+ & 1-\alpha & 0\\
0 & \alpha d_2^+ & 1-\alpha \\
1-\alpha & 0 & \alpha d_3^+
\end{bmatrix}
\]
with $d_1^+ , d_2^+, d_3^+ \ge 1 $ and $3 \le d_1^+ +d_2^+ +d_3^+ \le a$\\
For $\alpha\in [0, 1)$, $\mbox{Rank}(Q)=3$. Consequently $\mbox{Rank}(A_{\alpha}(D))$ is atleast three in this case, a contradiction.\\
Finally, If $D$ contains one of the induced subdigraphs $D_9,D_{10},D_{11}$, then the principal submatrix $R$ of $A_{\alpha}(D)$ corresponding to vertices in $X$ has the form
\[
{R}=
\begin{bmatrix} 
\alpha d_1^+ & 1-\alpha & 1-\alpha\\
0 & \alpha d_2^+ & 1-\alpha \\
0 & 0 & \alpha d_3^+
\end{bmatrix}
\]
with $d_1^+ \ge 2 , d_2^+ \ge 1, d_3^+ \ge 0 $ and $3 \le d_1^+ +d_2^+ +d_3^+ \le a$\\
If $\alpha=0$, then $\mbox{Rank}(R)=2$.
If $\alpha \ne 0$, then $\mbox{Rank}(R)=2$ or $3$ according as $d_3^+=0$ or $d_3^+ \ge 1$
Consequently, $\mbox{Rank}(A_{\alpha}(D))$ is atleast $2$ in this case, a contradiction.\\
Hence, $D$ cannot have two consecutive arcs $(v_1,v_2)$ and $(v_2,v_3)$. This proves the claim.\\
Now $D$ is either ${\overrightarrow{K}}_{r,s}$ or a proper subdigraph of ${\overrightarrow{K}}_{r,s}$ with no isolated vertex. Since $\alpha$-matrix of ${\overrightarrow{K}}_{r,s}$ has the form \[ A_{\alpha}({\overrightarrow{K}}_{r,s})=\begin{bmatrix} 
\alpha sI_{r} & (1-\alpha)J_{r \times s}\\
0_{s\times r}& 0_{s \times s}
\end{bmatrix}
\]

It is easy to see that any proper subdigraph of ${\overrightarrow{K}}_{r,s}$ with no isolated vertex has rank atleast two. Hence only possibility is that $D={\overrightarrow{K}}_{r,s}$. From Lemma $2.3$, we see $Rank(A_\alpha({\overrightarrow{K}}_{r,s}))=1$ if and only if $\alpha=0$, by noting that rank of a matrix equals to number of its non zero singular values.\qed

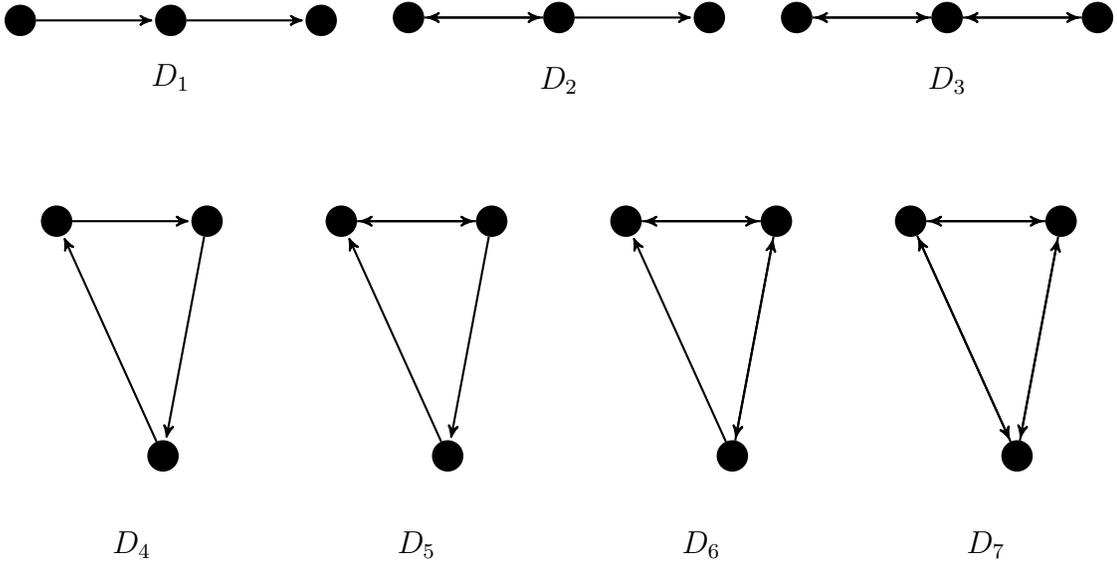
\begin{center}
\begin{minipage}{0.3\textwidth}
\centering
\begin{tikzpicture}[->,>=stealth',shorten >=1pt,auto,node distance=2cm,
                    thick,main node/.style={circle,draw,fill=black,minimum size=0.3cm}]

  \node[main node] (A) {};
  \node[main node] (B) [right of=A] {};
  \node[main node] (C) [right of=B] {};

  \path[every node/.style={font=\sffamily\small}]
    (A) edge node {} (B)
    (B) edge node {} (C);
\end{tikzpicture}

$D_1$
\end{minipage}%
\begin{minipage}{0.3\textwidth}
\centering
\begin{tikzpicture}[->,>=stealth',shorten >=1pt,auto,node distance=2cm,
                    thick,main node/.style={circle,draw,fill=black,minimum size=0.3cm}]

  \node[main node] (A) {};
  \node[main node] (B) [right of=A] {};
  \node[main node] (C) [right of=B] {};

  \path[every node/.style={font=\sffamily\small}]
    (A) edge node {} (B)
    (B) edge node {} (A)
    (B) edge node {} (C);
\end{tikzpicture}

$D_2$
\end{minipage}%
\begin{minipage}{0.3\textwidth}
\centering
\begin{tikzpicture}[->,>=stealth',shorten >=1pt,auto,node distance=2cm,
                    thick,main node/.style={circle,draw,fill=black,minimum size=0.3cm}]

  \node[main node] (A) {};
  \node[main node] (B) [right of=A] {};
  \node[main node] (C) [right of=B] {};

  \path[every node/.style={font=\sffamily\small}]
    (A) edge node {} (B)
    (B) edge node {} (A)
    (B) edge node {} (C)
    (C) edge node {} (B);
\end{tikzpicture}

$D_3$
\end{minipage}
\end{center}
\vskip 1mm
\begin{center}
\begin{minipage}{0.22\textwidth}
\centering
$$
\begin{tikzpicture}[->,>=stealth',shorten >=1pt,auto,node distance=2cm,
                    thick,main node/.style={circle,draw,fill=black,minimum size=0.1cm}]
  \node[main node] (A) {};
  \node[main node] (B) [right of=A] {};
  \node[main node] (C) [below right of=A, below=1.5cm] {};

  \path[every node/.style={font=\sffamily\small}]
    (A) edge node {} (B)
    (B) edge node {} (C)
    (C) edge node {} (A);
\end{tikzpicture}
$$
$D_4$
\end{minipage}%
\begin{minipage}{0.22\textwidth}
\centering
$$
\begin{tikzpicture}[->,>=stealth',shorten >=1pt,auto,node distance=2cm,
                    thick,main node/.style={circle,draw,fill=black,minimum size=0.3cm}]
  \node[main node] (A) {};
  \node[main node] (B) [right of=A] {};
  \node[main node] (C) [below right of=A, below=1.5cm] {};

  \path[every node/.style={font=\sffamily\small}]
    (A) edge node {} (B)
    (B) edge node {} (A)
    (B) edge node {} (C)
    (C) edge node {} (A);
\end{tikzpicture}
$$
$D_5$
\end{minipage}%
\begin{minipage}{0.22\textwidth}
\centering
$$
\begin{tikzpicture}[->,>=stealth',shorten >=1pt,auto,node distance=2cm,
                    thick,main node/.style={circle,draw,fill=black,minimum size=0.3cm}]
  \node[main node] (A) {};
  \node[main node] (B) [right of=A] {};
  \node[main node] (C) [below right of=A, below=1.5cm] {};

  \path[every node/.style={font=\sffamily\small}]
    (A) edge node {} (B)
    (B) edge node {} (A)
    (B) edge node {} (C)
    (C) edge node {} (B)
    (C) edge node {} (A);
\end{tikzpicture}
$$
$D_6$
\end{minipage}%
\begin{minipage}{0.22\textwidth}
\centering
$$
\begin{tikzpicture}[->,>=stealth',shorten >=1pt,auto,node distance=2cm,
                    thick,main node/.style={circle,draw,fill=black,minimum size=0.3cm}]
  \node[main node] (A) {};
  \node[main node] (B) [right of=A] {};
  \node[main node] (C) [below right of=A, below=1.5cm] {};

  \path[every node/.style={font=\sffamily\small}]
    (A) edge node {} (B)
    (B) edge node {} (A)
    (B) edge node {} (C)
    (C) edge node {} (B)
    (A) edge node {} (C)
    (C) edge node {} (A);
\end{tikzpicture}
$$
$D_7$
\end{minipage}
\end{center}
\vskip 3mm
\begin{center}
\begin{minipage}{0.23\textwidth}
\centering
$$
\begin{tikzpicture}[->,>=stealth',shorten >=1pt,auto,node distance=2cm,
                    thick,main node/.style={circle,draw,fill=black,minimum size=0.3cm}]
  \node[main node] (A) {};
  \node[main node] (B) [right of=A] {};
  \node[main node] (C) [below right of=A, below=1.5cm] {};

  \path[every node/.style={font=\sffamily\small}]
    (A) edge node {} (B)
    (B) edge node {} (C)
    (A) edge node {} (C);
\end{tikzpicture}
$$
$D_8$
\end{minipage}%
\begin{minipage}{0.23\textwidth}
\centering
$$
\begin{tikzpicture}[->,>=stealth',shorten >=1pt,auto,node distance=2cm,
                    thick,main node/.style={circle,draw,fill=black,minimum size=0.3cm}]
  \node[main node] (A) {};
  \node[main node] (B) [right of=A] {};
  \node[main node] (C) [below right of=A, below=1.5cm] {};

  \path[every node/.style={font=\sffamily\small}]
    (A) edge node {} (B)
    (B) edge node {} (C)
    (C) edge node {} (B)
    (A) edge node {} (C);
\end{tikzpicture}
$$
$D_9$
\end{minipage}%
\begin{minipage}{0.23\textwidth}
\centering
$$
\begin{tikzpicture}[->,>=stealth',shorten >=1pt,auto,node distance=2cm,
                    thick,main node/.style={circle,draw,fill=black,minimum size=0.3cm}]
  \node[main node] (A) {};
  \node[main node] (B) [right of=A] {};
  \node[main node] (C) [below right of=A, below=1.5cm] {};

  \path[every node/.style={font=\sffamily\small}]
    (A) edge node {} (B)
    (B) edge node {} (A)
    (B) edge node {} (C)
    (A) edge node {} (C);
\end{tikzpicture}
$$
$D_{10}$
\end{minipage}%
\begin{minipage}{0.23\textwidth}
\centering
$$
\begin{tikzpicture}[->,>=stealth',shorten >=1pt,auto,node distance=2cm,
                    thick,main node/.style={circle,draw,fill=black,minimum size=0.3cm}]
  \node[main node] (A) {};
  \node[main node] (B) [right of=A] {};
  \node[main node] (C) [below right of=A, below=1.5cm] {};

  \path[every node/.style={font=\sffamily\small}]
    (A) edge node {} (B)
    (B) edge node {} (A)
    (B) edge node {} (C)
    (C) edge node {} (B)
    (A) edge node {} (C);
\end{tikzpicture}
$$
$D_{11}$
\end{minipage}
\end{center}
\vskip 1mm
$$Fig. 1$$

Recall a digraph is said to be discrete if it has no arcs. Using the definition and singular value decomposition, the following result holds 
\begin{lemma}\label{2.6}
 D is a discrete digraph of order $n$ if and only if $\alpha$-singular values of D are $0^{[n]}$.
\end{lemma}
Let $D_1=(\mathcal{V}_1,\mathcal{A}_1), D_2=(\mathcal{V}_2,\mathcal{A}_2),\dots, D_k=(\mathcal{V}_k,\mathcal{A}_k)$ be $k$ digraphs. Then their direct sum $D=\oplus_{i=1}^{k}D_i$ is digraph with vertex set $\mathcal{V}=\cup_{i=1}^{k}\mathcal{V}_i$ and arc set $\mathcal{A}=\cup_{i=1}^{k}\mathcal{A}_i$. The following result is easy to prove.
\begin{lemma}\label{2.7}
 If $D$ is a direct sum of $D_1,D_2,...,D_k$ digraphs, then $$ \|{D_{\alpha}} \|_* =\sum_{i=1}^{k} \|{(D_i)_{\alpha}} \|_* $$.
\end{lemma}
\begin{lemma}\label{2.8}
 If $D$ is a k-regular digraph, then each row sum of $A_{\alpha}{A_{\alpha}}^{T}$ equals $k^2$.
 \end{lemma}
 {\bf Proof.} For a $k$-regular digraph 
 \begin{align*}
 A_{\alpha}{A_{\alpha}}^{T}&=[\alpha kI_n+(1-\alpha)A][\alpha kI_n+(1-\alpha)A^T]\\
 &=\alpha^2k^2I_n+k\alpha(1-\alpha)(A+A^T)+(1-\alpha)^2AA^T.
 \end{align*}
 We first compute the row sum of $A_{\alpha}{A_{\alpha}}^{T}$. Note that the ith row sum of $AA^T$ is equal to $\sum_{j=1}^{n} \langle R_i, R_j \rangle  $. In each row $R_i$ of $A$, there are $k$ non-zero entries equal to $1$. Assume the ones are at the positions $i_1,i_2,\cdots i_k$.\\
For each $i_l$ where $1\le l \le k$ as  $d^{-}(v_{i_l})=k$, there are $k-1$ ones above and below entry $[A]_{i,i_l}$.\\
Therefore, \\
 \begin{align*}
 \sum_{j=1}^{n} \langle R_i, R_j \rangle  &=\langle R_i, R_i \rangle  +\sum_{j \ne i} \langle R_i, R_j \rangle  \\
&=k+k(k-1)\\
 &=k^2
 \end{align*}
 Now, \begin{align*} \sum_{j=1}^{n}[A_{\alpha}{A_{\alpha}}^{T}]_{ij}&=\alpha^2k^2+k\alpha(1-\alpha)(2k)+(1-\alpha)^2k^2\\
 &=\alpha^{2}k^{2}+2k^{2}\alpha-2k^{2}\alpha^{2}+k^{2}+k^{2}\alpha^{2}-2k^{2}\alpha\\
 &=k^{2}. 
  \end{align*}\qed
  
We recall that if $\sigma_{1\alpha},\sigma_{2\alpha},\dots,\sigma_{n\alpha}$ are singular values of $\alpha$ matrix $A_{\alpha}(D)$ of a digraph $D$, then $\sigma_{1\alpha}$ is known as the $\alpha$-spectral norm of $D$. Cruz, Giraldo and Rada \cite{cgr} obtained a lower bound for the spectral norm of a digraph. We next determine a lower bound for $\alpha$ spectral norm, independent of $\alpha$.
\begin{lemma}\label{2.9}
 If $D$ is a digraph with $n\ge2$ vertices and $a$ arcs, then $\sigma_{1\alpha}(D)\ge \frac{a}{n} $ with equality if and only if $D$ is $\frac{a}{n}$-regular digraph.
 \end{lemma}
  {\bf Proof.} We recall $\|.\|_2$ denote usual Euclidean norm. Let ${\bf{e}}\in {\mathbb{R}}^n$ be the column vector with all entries $1$. As $A_{\alpha}{A_{\alpha}}^T$ is real symmetric matrix, using Rayleigh quotient and Cauchy Schwarz inequality, we see  
   \begin{align*}
(\sigma_{1\alpha}(D))^2&=\max_{x \ne 0}\frac{x^{T}A_{\alpha}A_{\alpha}^{T}x}{x^{T}x}
 =\max_{x \ne 0}\frac{x^{T}A_{\alpha}^{T}A_{\alpha}x}{x^{T}x}\\
 &=\max_{x \ne 0}\frac{\|A_{\alpha}x\|_{2}^{2}}{\|x\|_{2}^{2}}\ge \frac{\|A_{\alpha}e\|_{2}^{2}}{\|e\|_{2}^{2}}\\
 &=\frac{\sum_{j=1}^{n}(d_{j}^{+})^2}{n}\ge \frac{ (\sum_{j=1}^{n}d_{j}^+)^{2}}{n^2}=\frac{a^{2}}{n^{2}}.
   \end{align*}
 Equality holds if and only if $d_j^{+}=\mbox{constant}=d^{+}$, say for all  $j=1,2,\cdots, n$ and ${\bf{e}}$ is an eigenvector of $A_{\alpha}^T A_{\alpha}$ i.e., row sums of $A_{\alpha}^T A_{\alpha}$ are constant.\\
As $D$ is outregular (i.e., all outdegrees of $D$ are equal), therefore $A_{\alpha}{\bf{e}}=d^{+}{\bf{e}}$.\\
Now row sums of $A_{\alpha}^T A_{\alpha}$ are constant gives $A_{\alpha}^T[d^{+}{\bf{e}}]=c {\bf{e}}$,  where $c$ is some constant.\\
or $\alpha (d^{+})^2 {\bf{e}}+(1-\alpha)d^{+}{\bf{d}}=c {\bf{e}}$, where ${\bf{d}}=[d_1^{-},d_2^{-},\dots,d_n^{-}]\in\mathbb{R}^n$ is indegree column vector of $D$.\\
In particular for $i\neq j$, we have\\
$$\alpha (d^{+})^2+(1-\alpha)d^{+}d_i^{-}=\alpha (d^{+})^2+(1-\alpha)d^{+}d_j^{-},$$
which upon simplification gives $d_i^{-}=d_j^{-}$ for all $i,j=1,2,\dots,n$ and $i\neq j$. This implies $D$ is indegree regular (i.e., all indegrees of $D$ are equal) as well. 
 Since $$\sum_{j=1}^{n}d_{j}^{+}=\sum_{j=1}^{n}d_{j}^{-},$$ we see that $d_j^{+}=d_j^{-}=\frac{a}{n}$ for all $j=1,2,\dots,n$, which in turn implies $D$ is $\frac{a}{n}$-regular digraph.\qed
We next obtain a lower bound for the $\alpha$ trace norm of digraphs.
\begin{theorem}\label{2.10}
 If $D$ is a digraph with $n\ge 2$ vertices and $a$ arcs and $A_{\alpha}=\alpha {\Delta^{+}}+(1-\alpha)A$ be the $\alpha$ matrix of $D$ and $\sigma_{1{\alpha}}\ge \sigma_{2{\alpha}} \ge\dots \ge \sigma_{n{\alpha}}\ge 0$ be the $\alpha$-singular values of $D$. Then 
 \begin{align} \|{D_{\alpha}} \|_*\ge \sqrt{(1-\alpha)^2a+\alpha^2\sum_{i=1}^{n}(d_i^{+})^2+n(n-1)|detA_{\alpha}|^{\frac{2}{n}}}
 \end{align}
 with equality if and only if\\
 (a) $D$ is a Discrete digraph~~ (b) $\alpha=0$ and $D$ is the direct sum of directed cycles.\\
 (c) $\alpha$=0 and $D=\overrightarrow{K}_{r,s}+$ Possibly some isolated vertices or $\alpha=\frac{1}{2}$ and  $D=\overleftrightarrow{K_{2}}.$\\
 \end{theorem}
 {\bf Proof.} Let $\sigma_{1{\alpha}}\ge \sigma_{2{\alpha}} \ge,\dots \ge \sigma_{n{\alpha}}\ge 0$ be the singular values of $A_{\alpha}=\alpha {\Delta^{+}}+(1-\alpha)A$. By AM-GM inequality 
 \begin{align*}
 \|{D_{\alpha}} \|_*^2&=(\sum_{i=1}^{n}\sigma_{i{\alpha}})^2\\
 &=\sum_{i=1}^{n}(\sigma_{i{\alpha}})^2+2\sum_{1\le i <j\le n}\sigma_{i{\alpha}}\sigma_{j{\alpha}}
 \\&=(1-\alpha)^2\sum_{i=1}^{n}d_i^{+}+\alpha^2\sum_{i=1}^{n}(d_i^{+})^2+\frac{2n(n-1)}{n(n-1)}\sum_{1\le i <j\le n}\sigma_{i{\alpha}}\sigma_{j{\alpha}}
\end{align*}
 
\begin{align*}
 &\ge(1-\alpha)^2a+\alpha^2\sum_{i=1}^{n}(d_i^{+})^2+n(n-1)(\prod_{i=1}^n(\sigma_{i{\alpha}})^{n-1})^\frac{2}{n(n-1)}\\
  &=(1-\alpha)^2a+\alpha^2\sum_{i=1}^{n}(d_i^{+})^2+n(n-1)|detA_{\alpha}|^{\frac{2}{n}}\\
   \|{D_{\alpha}} \|_*&\ge \sqrt{(1-\alpha)^2a+\alpha^2\sum_{i=1}^{n}(d_i^{+})^2+n(n-1)|detA_{\alpha}|^{\frac{2}{n}}}
  \end{align*}
   This proves (2.1).\\
   Assume equality holds in (2.1), then equality holds in AM-GM inequality which gives 
   $$\sigma_{i{\alpha}}\sigma_{j{\alpha}}=c $$ for all $i,j$. Three cases depending on $\alpha$- singular values arise here\\
   {\bf Case (1)}. If $\sigma_{1{\alpha}}=0$, then 
   $\sigma_{1{\alpha}}=\sigma_{2{\alpha}}=\dots \sigma_{3{\alpha}}=0$. This implies $D$ is a Discrete digraph, by Lemma $2.6$. \\
   {\bf Case (2)}. If $\sigma_{1{\alpha}}>0,\sigma_{2{\alpha}}>0$. Then $\sigma_{1{\alpha}}=\sigma_{2{\alpha}}= \dots =\sigma_{n{\alpha}}$. By singular value decomposition, there exists a real orthogonal matrix $U$ such that \\
      \begin{align*}
      A_\alpha&=\sigma_{1{\alpha}}{U}\\
      \mbox {or}~~~ [A_\alpha]_{ij}&=\sigma_{1{\alpha}}[{U}]_{i,j}
      \end{align*}
      
      \begin{align*}
     \mbox {or}~~~    \sum_{j=1}^n[A_\alpha]_{ij}^2&= {(\sigma_{1{\alpha}})}^2\sum_{j=1}^n{[{U}]_{ij}}^2
     =(\sigma_{1{\alpha}})^2\\
     \mbox{or} ~~~~~~~\alpha^2(d_i^{+})^2+(1-\alpha)^2d_i^{+}&=(\sigma_{{1\alpha}})^2~~
\mbox {for all}~~~ i = 1, 2, 3, \dots, n
\end{align*}
Again \begin{align*}\sum_{i=1}^n[A_\alpha]_{ij}^2&= {(\sigma_{1{\alpha}})}^2\sum_{i=1}^n{[{U}]_{ij}}^2\\
\mbox{or} ~~~~~~~\alpha^2(d_j^{+})^2+(1-\alpha)^2d_j^{-}&=(\sigma_{1{\alpha}})^2~~~~~~
\mbox{for all}~~ j = 1, 2, 3, \dots, n  
\end{align*}

In particular $$\alpha^2 (d_i^{+})^{2}+(1-\alpha)^2 d_i^{+}
=\alpha^2 (d_i^{+})^{2}+(1-\alpha)^2d_i^{-}$$
which gives $d_i^+=d_i^-=d$, say~~for all ~~$i=1,2, \dots,n.$
This implies $D$ is $d$-regular.\\
Using Lemma $2.9$, we see, \begin{align*}\sum_{i=1}^n (\sigma_{i{\alpha}})^2 &=(1-\alpha)^2a+\alpha^2\sum_{i=1}^n(d_i^{+})^2\\
 n(\sigma_{i{\alpha}})^2&=(1-\alpha)^2nd+\alpha^2nd^2\\
nd^2&=(1-\alpha)^2nd+\alpha^2nd^2\\
(1-{\alpha}^2)nd^2&=(1-\alpha)^2nd\\
d&=\frac{1-\alpha}{1+\alpha}
\end{align*}
which gives $$ \alpha=0~~~~and ~~d=1$$
So, $D$ is direct sum of directed cycles and $\alpha=0$.\\
{\bf Case (3)}.
   If $\sigma_{1{\alpha}}>0,\sigma_{2{\alpha}}=0$, then 
    $\sigma_{1{\alpha}}>0,\sigma_{2{\alpha}}=\sigma_{3{\alpha}}=\dots \sigma_{n{\alpha}}=0$.\\
By singular value decomposition, there exists real orthogonal matrices $U=(u_{ij})_{n\times n}$ and $V=(v_{ij})_{n\times n}$ such that
\begin{align*}
A_{\alpha}={U}S {V^{T}},
\end{align*}
  where  $ S=diag( \sigma_{1\alpha},0,0,\ldots,0)
   =diag(\sqrt{c},0,0,\ldots,0),$ where $c=(1-\alpha)^2a+\alpha^2\sum_{i=1}^n(d_i^{+})^2$\\
   In this case 
   $A_{\alpha}=\sqrt{c}{u}{v}^T,$
 where $u=[u_{11},~~u_{21},~~\dots,u_{n1}]^T$ and $v=[v_{11},~~ v_{21},~~\dots,v_{n1}]^T$.\\ 
This implies $A_\alpha$  is a rank $1$ matrix. By theorem $2.5$, we see that $\alpha=0$ and $D=\overrightarrow{K}_{r,s}$ + Possibly some isolated vertices or $\alpha=\frac{1}{2}$ and  $D=\overleftrightarrow{K_{2}}$. \\
\indent Conversely, if $D$ is discrete digraph, then both sides of $(2.1)$ are equal to zero. If $\alpha=0$ and $D$ is direct sum of cycles, then both sides of $(2.1)$ are equal to $a$, the number of arcs in $D$. If $D=\overrightarrow{K}_{r,s}+$ Possibly some isolated vertices, the both sides of $(2.1)$ are equal to $\sqrt{rs}$. If $\alpha=\frac{1}{2}$ and  $D=\overleftrightarrow{K_{2}}$, then both sides of $(2.1)$ are equal to $1$. This completes the proof.\qed

Following is an immediate consequence of Theorem $2.10$.\\
\begin{corollary}\label{2.11}
If $D$ is a digraph with $n\ge 2$ vertices and $a$ arcs then $$\|{D_\alpha} \|_*\ge \sqrt{(1-\alpha)^2a+\alpha^2\sum_{i=1}^n(d_i^{+})^2},$$
with equality if and only if $D$ satisfies one of the following\\
~(a)~$D$ is a discrete digraph.\\ (b)~$\alpha=0$ and $D=\overrightarrow{K}_{r,s}+$ Possibly some isolated vertices or $\alpha=\frac{1}{2}$ and $D=\overleftrightarrow{K_2}$. 
\end{corollary}
Given a tree $T$ on $n$, vertices, let $\mathcal{T}(n)$ denote the set of oriented trees with underlying tree $T$. The next result determines oriented trees having minimum $\alpha$-trace norm in $\mathcal{T}(n)$.
\begin{corollary}\label{2.12}
 If $T\in \mathcal{T}(n)$, then \begin{align} \|{T_{\alpha}} \|_* \ge \sqrt{(1-\alpha)^2(n-1)+\alpha^2\sum_{i=1}^{n}(d_i^{+})^2},\end{align} 
 with equality if and only if $\alpha=0$ and $T={\overrightarrow{
 K}_{1,n-1}} ~\mbox{or}~{\overrightarrow{
 K}_{n-1,1}} $.
 \end{corollary}
  {\bf Proof.}
  We note that for a digraph $D$
  $$ \|{D_{\alpha}} \|_* \ge \sqrt{(1-\alpha)^2a+\alpha^2\sum_{i=1}^{n}(d_i^{+})^2}$$ 
  with equality if and only if $\alpha=0$ and $D={\overrightarrow{K}_{r,s}}$ + Possibly some isolated vertices or $\alpha=\frac{1}{2}$ and $D=\overleftrightarrow{K_2}$.\\
  For an oriented tree $T\in \mathcal{T}(n)$, we see 
   $$ \|{T_{\alpha}} \|_* \ge \sqrt{(1-\alpha)^2(n-1)+\alpha^2\sum_{i=1}^{n}(d_i^{+})^2}.$$ 
Also, since $T$ has no isolated vertices and no cycles and ${\overrightarrow{K}_{r,s}}$ is a tree if and only if either $r=1$ and $s=n-1$ or $r=n-1$ and $s=1$, the result follows.\qed
\section{\bf Upper bounds for trace norm of $A_{\alpha}$ matrix of digraphs}\label{sec3}
The well known McClelland's upper bound [Theorem 5.1, \cite{lsg}] states that for a graph $G$ with $p$ vertices and $q$ edges
 \begin{eqnarray}
E(G)\le \sqrt{2pq}.\end{eqnarray}
Moreover, equality holds in $(3.1)$ if and only if $G$ is a discrete graph or $G=(\frac{n}{2})K_2$.\\
This upper bound has been extended to digraph energy by Rada \cite{r1}. We next find McClelland type upper bound for trace norm of $A_{\alpha}$ matrix of a digraph. Recall that for nonnegative real numbers $x_1,x_2,\dots,x_n$, their variance $var(x_1,x_2,\dots,x_n)\ge 0$ with equality if and only if $x_1=x_2=\dots=x_n$.
\begin{theorem} \label{3.1} Let $D$ be a digraph with $n$ vertices, $a$ arcs and let $(d_1^+,d_2^+,\dots,d_n^+)$ be out degrees of vertices. Then 
\begin{align}
 \|{D} \|_*  \le \sqrt{n[(1-\alpha)^2a+\alpha^2\sum_{i=1}^{n}(d_i^{+})^2]},
 \end{align}
 with equality if and only if $D$ satisfies one of the following\\
 (a) $D$ is a Discrete digraph.\\
 (b) $\alpha=0$ and $D$ is a direct sum of directed cycles.
\end{theorem}
{\bf Proof.}
Let $\sigma_{1\alpha}\ge \sigma_{2\alpha} \ge,\dots \ge \sigma_{n\alpha}\ge 0$ be the singular values of $A_{\alpha}=\alpha \Delta^{+}+(1-\alpha)A$\\
We know that 
\begin{align*}
0&\le var[{\sigma_{1\alpha},\sigma_{2\alpha},\sigma_{3\alpha},\dots,\sigma_{n\alpha}}]\\
&=\frac{1}{n}\sum_{i=1}^n (\sigma_{i\alpha})^2-(\frac{\sum_{i=1}^n\sigma_{i\alpha}}{n})^2 \\
&=\frac{1}{n}\sum_{i=1}^n (\sigma_{i\alpha})^2-(\frac{\sum_{i=1}^n\sigma_{i\alpha}}{n})^2 \\
&=\frac{1}{n}[(1-\alpha)^2a+ \alpha^2\sum_{i=1}^n (d_i^+)^2]-\frac{\|{D_{\alpha}} \|_* ^2}{n^2} \\
\mbox{or}~~~\frac{\|{D_{\alpha}} \|_* ^2}{n^2} &\le \frac{1}{n}[(1-\alpha)^2a+ \alpha^2\sum_{i=1}^n (d_i^+)^2]
\end{align*}

\begin{align*}
\mbox{or}~~~\frac{\|{D_{\alpha}} \|_* ^2}{n}& \le [(1-\alpha)^2a+ \alpha^2\sum_{i=1}^n (d_i^+)^2]\\
\mbox{or}~~~{\|{D_{\alpha}} \|_* ^2}& \le n[(1-\alpha)^2a+ \alpha^2\sum_{i=1}^n (d_i^+)^2]\\
\mbox{or}~~~{\|{D_{\alpha}} \|_* }& \le \sqrt{n[(1-\alpha)^2a+ \alpha^2\sum_{i=1}^n (d_i^+)^2]}
\end{align*}
This proves the inequality. The equality holds in $(3.2)$ if and only if $\sigma_{1\alpha}= \sigma_{2\alpha} =\dots = \sigma_{n\alpha}.$ Proceeding as in the equality case of Theorem $2.10$, we see equality holds if and only if $D$ is a discrete digraph or $\alpha=0$ and $D$ is direct sum of cycles. \qed
\begin{remark}\label{3.2} If $D=\overleftrightarrow{G}$ is a symmetric digraph on $n$ vertices, then for $\alpha=0$, upper bound $(3.2)$ reduces to McClelland's upper bound for graphs. Therefore, Theorem $3.1$, is an extension of McClelland's inequality for graphs.
\end{remark}
 
A square $(0,1)$ matrix $M$ of order $n$ is the incidence matrix of a symmetric $(n,k,\lambda)$-BIBD if and only if$$MM^T=\lambda J+(k-\lambda)I,$$ 
where $J$ is matrix of all ones and $I$ stands for identity matrix. If adjacency matrix of a digraph $D$ satisfies this condition, we say $D$ is a symmetric $(n,k,\lambda)$-BIBD.\\

A $k$ regular graph $G$ of order $n$ is said to be strongly regular graph $\mbox{(SRG)}$ with parameters $(n,k,\lambda,\mu)$ if any two adjacent vertices of $G$ have $\lambda$ common neighbours and any two non adjacent vertices of $G$ have $\mu$ common neighbours. For example Shirkhande graph is a strongly regular graph $\mbox{SRG}~(16,6,2,2)$. Since its two parameters $\lambda\mbox~{and}~\mu$ are equal, its adjacency matrix clearly satisfies condition for being a symmetric BIBD.\\ 
 
We next obtain Koolen and Moulton type upper bound \cite{km} for the $\alpha$ trace norm of digraphs. We adapt the idea from \cite{n1} and discuss the equality case. Here by $d^{+}_{\mbox{max}}$, we denote the maximum outdegree of digraph $D$.
\begin{theorem}\label{3.3}
 Let $D$ be a digraph with $n\ge 2$ vertices and $a\ge n\beta$ arcs, \\where $\beta=\max{(1-\alpha,\alpha~ d^{+}_{\mbox{max}})}$. Then
 \begin{align}
 \|D_{\alpha}\|_{*}\le \frac{a}{n}+\sqrt{(n-1)[(1-\alpha)^2a+\alpha^2\sum_{i=1}^{n}(d_i^{+})^2-\frac{a^2}{n^2}]}
 \end{align}
 with equality if and only if $D$ satisfies one of the following\\
 (a) $\alpha=0$ and $D$ is direct sum of cycles.
 (b) $D$ is $\frac{a}{n}$-regular with two distinct nonnegative singular values and these are $\sigma_{1\alpha}=\frac{a}{n},\sigma_{2\alpha}=\sigma_{3\alpha}=\dots=\sigma_{n\alpha}=\sigma=\sqrt{\frac{\alpha^2\sum_{i=1}^{n}(d_i^+)^2+(1-\alpha)^2 a-\frac{a^2}{n^2}}{n-1}}.$
 \end{theorem}
  {\bf Proof.} Let $\sigma_{1\alpha},\sigma_{2\alpha},\dots,\sigma_{n\alpha}$ be the $\alpha$-singular values of $D$. Then by Cauchy Schwarz Inequality, we have 
  \begin{align*}
  \|D_{\alpha}\|_{*}& =\sum_{i=1}^{n}\sigma_{i\alpha}
=\sigma_{1\alpha}+\sum_{i=2}^{n}\sigma_{i\alpha}\\
&\le \sigma_{1\alpha}+\sqrt{(n-1)\sum_{i=2}^{n}(\sigma_{i\alpha})^{2}}
\\&=\sigma_{1\alpha}+\sqrt{(n-1)[(1-\alpha)^2a+\alpha^{2}\sum_{i=2}^{n}(d_{i}^+)^2-(\sigma_{1\alpha})^{2}]}
\end{align*}
Consider the function
$$f(x)=x+\sqrt{(n-1)[(1-\alpha)^2a+\alpha^{2}\sum_{i=2}^{n}(d_{i}^+)^2-x^{2}]}$$
on $\left[0, \sqrt{(n-1)\left[(1-\alpha)^2 a + \alpha^{2} \sum_{i=1}^{n} (d_{i}^+)^2\right]}~ \right]
$. It is easy to verify that $f$ is strictly decreasing on $$\left[\sqrt{\frac{(1-\alpha)^2a+\alpha^{2}\sum_{i=1}^{n}(d_{i}^+)^2}{n}},\sqrt{(1-\alpha)^2a+\alpha^{2}\sum_{i=1}^{n}(d_{i}^{+})^2)} ~\right].$$
Also, with $a\ge n\beta$, we have \\
  \begin{align*}
(1-\alpha)^{2}a+\alpha^{2}\sum_{i=1}^{n}(d_{i}^{+})^2&=Trace(A_{\alpha}A_{\alpha}^{*})\\&= \sum_{i=1}^{n}\sum_{j=1}^{n}[A_{\alpha}]_{ij}^{2}\\&\le \beta \sum_{i=1}^{n}\sum_{j=1}^{n}[A_{\alpha}]_{ij}=\beta a
\end{align*}

\begin{align*}
\mbox{or} ~~~~(1-\alpha)^{2}a+\alpha^{2}\sum_{i=1}^{n}(d_{i}^{+})^2 &\le \beta a\le \frac{a^2}{n}\\
\mbox{or}~~~~\frac{(1-\alpha)^{2}a+\alpha^{2}\sum_{i=1}^{n}(d_{i}^{+})^2}{n}&\le  \frac{a^{2}}{n^{2}}\\
\mbox{or}~~~~\sqrt{\frac{(1-\alpha)^{2}a+\alpha^{2}\sum_{i=1}^{n}(d_{i}^{+})^2}{n}}&\le  \frac{a}{n}\le \sigma_{1\alpha} ~~\text[{By~~ Lemma~~ 2.9}]
  \end{align*}
Therefore, $$\|D_{\alpha}\|_{*}\le f(\sigma_{1\alpha})\le f(\frac{a}{n})$$ which proves bound (3.3).\\
Moreover, equality holds if and only if $D$ has at most two singular values $\sigma_{1\alpha},\sigma_{2\alpha}=\sigma_{3\alpha}=\dots=\sigma_{n\alpha}$ and $D$ is $\frac{a}{n}$-regular so $\sigma_{1\alpha}=\frac{a}{n}$.\\
In case $D$ has only one singular value then $\sigma_{1\alpha}=\sigma_{2\alpha}=\dots=\sigma_{n\alpha}$. Then since $a\ge n\beta$, proceeding as in equality case of Theorem $2.10$, $\alpha=0$ and $D$ is direct sum of cycles. In view of Theorem $2.5$, we see that if $\sigma_{i\alpha}=0$, for $i=2,3,\dots,n$, then $n=2$, $\alpha=\frac{1}{2}$ and $D=\overleftrightarrow{K_2}$. Thus, we conclude that if $D$ has two distinct positive singular values, then these are $\sigma_{1\alpha}=\frac{a}{n},\sigma_{2\alpha}=\sigma_{3\alpha}=\dots,\sigma_{n\alpha}=\sigma=\sqrt{\frac{\alpha^2\sum_{i=1}^{n}(d_i^+)^2+(1-\alpha)^2 a-\frac{a^2}{n^2}}{n-1}}$. Converse part can be easily verified.\qed
\begin{remark}\label{3.4}
At this moment, we are not able to determine digraphs for which $A_{\alpha}$ matrix has two distinct positive singular values of the form stated in Theorem $3.3$. It remains a future problem to determine such digraphs. Here we note that for $\alpha=0$, if a digraph $D$ attain upper bound given in Theorem $3.3$, then $D$ must be a symmetric $(n,\frac{a}{n}, \frac{a(a-n)}{n^2(n-1})$-BIBD. For example see [Theorem 1.1, \cite{ar}]. For other values of $\alpha$, these digraphs need not have the desired property. For example take symmetric digraph associated with Shirkhande graph which is a strongly regular graph $\mbox{SRG}~(16,6,2,2)$. Its eigenvalues are $6, 2^{[6]},-2^{[9]}$. The singular values of its $\alpha$ matrix are $6, (2+4\alpha)^{[6]},|8\alpha-2|^{[9]}$. We observe that for $\alpha=0$, it has two distinct singular values. It is easy to see that for $\alpha\in[0, 1)$, the singular values of $\overleftrightarrow{K}_n$ are $n-1,|n\alpha-1|^{[n-1]}$. So, for $\alpha$ matrix of $\overleftrightarrow{K}_n$, equality holds in upper bound give in Theorem $3.3$.
\end{remark}
\noindent{\bf Acknowledgements.} The research of Mushtaq A. Bhat is supported by SERB-DST grant with File No. MTR/2023/000201. The research of Peer Abdul Manan is supported by CSIR, New Delhi, India with CSIR-HRDG Ref. No: Jan-Feb/06/21(i)EU-V. This research is also supported by NBHM project number NBHM/02011/20/2022.

\end{document}